\documentclass[12pt,draft]{amsart}
\usepackage{amsmath,amsthm,latexsym,amscd,amsbsy,amssymb,fleqn,leqno,pb-diagram}
\chardef\bslash=`\\ 

\makeatletter
\def\verbatim{\interlinepenalty\@M \@verbatim
  \leftskip\@totalleftmargin\advance\leftskip2pc
  \frenchspacing\@vobeyspaces \@xverbatim}
\makeatother
\makeatletter
  \def\dgt@k{\dg@DX=-3 \dg@DY=2 \dg@SIZE=3}
\makeatother

\makeatletter
  \def\dgt@kk{\dg@DX=3 \dg@DY=-1 \dg@SIZE=3}%
\makeatother

\theoremstyle{plain}
\newtheorem{thm}{Theorem}[section]
\newtheorem{cor}[thm]{Corollary}
\newtheorem{lem}[thm]{Lemma}
\newtheorem{pro}[thm]{Proposition}

\theoremstyle{definition}

\newtheorem*{remark}{Remark}

\numberwithin{equation}{section}

\newcounter{rmnum}

\def\symbolnote#1#2{\let\thefootn=\thefootnote%
\renewcommand{\thefootnote}{\fnsymbol{footnote}}%
\footnotemark[#1]%
\footnotetext[#1]{#2}%
\let\thefootnote=\thefootn
}

\newfont{\bbb}{msbm10 scaled \magstep1}
\newfont{\bbc}{msbm8 scaled \magstep0}

\newcommand{\R}{\mbox{\bbb R}}

\newcommand{\uin}{\mbox{\bbb I}}
\newcommand{\e}{\mbox{\rm e-dim}}

\newcommand{\St}{\hbox{\rm St}}

\newcommand{\mesh}{\hbox{\rm mesh}}

\newcommand{\cone}{\hbox{\rm Cone}}

\begin{document}


\title{Extension dimension and quasi-finite $CW$-complexes}

\author{Alex Karasev}
\address{Department of Computer Science and Mathematics,
Nipissing University,
100 College Drive, P.O. Box 5002, North Bay, ON, P1B 8L7, Canada}
\email{alexandk@nipissingu.ca}
\thanks{The first author was partially supported by his NSERC grant.}

\author{Vesko Valov}
\address{Department of Computer Science and Mathematics, Nipissing University,
100 College Drive, P.O. Box 5002, North Bay, ON, P1B 8L7, Canada}
\email{veskov@nipissingu.ca}
\thanks{The second author was partially supported by his NSERC grant.}

\keywords{Extension dimension, cohomological dimension, absolute
extensor, universal space, quasi-finite complex}
\subjclass{Primary: 54F45; Secondary: 55M10, 54C65.}


\begin{abstract}
We extend the definition of quasi-finite complexes by considering
not necessarily countable complexes. We provide a characterization
of quasi-finite complexes in terms of $L$-invertible maps and
dimensional properties of compactifications. Several results
related to the class of quasi-finite complexes are established,
such as completion of metrizable spaces, existence of universal
spaces and a version of the factorization theorem. Further, we
extend the definition of $UV(L)$-spaces on non-compact case and
show that some properties of $UV(n)$-spaces and $UV(n)$-maps
remain valid, respectively,  for $UV(L)$-spaces and $UV(L)$-maps.
\end{abstract}

\maketitle

\markboth{A.~Karasev and V.~Valov}{Quasi-finite $CW$-complexes}

\section{Introduction}
Extension theory introduced by Dranishnikov \cite{d,dd} unifies
the covering dimension and the cohomological dimension. There are
two classes of maps which play an important role in extension
theory. For a given complex $L$, theses are $L$-invertible and
$L$-soft maps. It should be mentioned that universal spaces in
dimension $L$ as well as absolute extensors in dimension $L$ are
obtained as preimages of Hilbert cube or Hilbert space under maps
from the above classes \cite{ch1}. For a countable complex $L$,
existence of $L$-invertible mapping of certain $L$-dimensional
compactum onto the Hilbert cube is closely connected with the
dimensional properties of compactifications of spaces with
extension dimension not grater than $L$ \cite{ch99}. It turned out
that the existence of such $L$-invertible mappings can be
characterized in terms of ``extensional" properties of a complex.
This inspired the concept of quasi-finite countable complexes
\cite{k:03}.

In the present paper we extend the definition of quasi-finite
complexes by considering not necessarily countable complexes. We
also provide a characterization of quasi-finite complexes in terms
of $L$-invertible maps and dimensional properties of
compactifications. Another interesting observation consists in the
fact that many results established for finite or countable
complexes remain valid for quasi-finite complexes. In particular,
quasi-finite complexes possess the $L$-soft map property and
every metrizable space of extension dimension $\leq L$ has a
completion with the same extensional dimension.  We also
prove a version of the factorization theorem, and construct
universal spaces. Finally, in case $L$ being quasi-finite it is
possible to define $UV(L)$-property for non-compact spaces. We
show that this property does not depend on the embedding of a
space into absolute neighborhood extensor in dimension $L$ and
obtain some results about $UV(L)$-maps and $UV(L)$-spaces which
were  known for $UV(n)$-maps and $UV(n)$-spaces, respectively.


\section{Quasi-finite $CW$-complexes}

Everywhere in this paper we assume that spaces are Tychonov and
maps are continuous. Let $X$ and $Y$  be two spaces, $A\subset X$
and $g\colon A\to Y$ a map. We write $Y\in ANE(g,A,X)$ if $g$ has
a continuous extension $\bar{g}\colon U\to Y$, where $U$ is a
neighborhood of $A$ in $X$ which has the following property: there
exists a function $h\colon X\to [0,1]$ such that $h^{-1}((0,1])=U$
and $h(A)=1$. If, in the above definition, $U=X$, we write $Y\in
AE(g,A,X)$. Let us note that, by \cite[Lemma 2.8]{dy1}, $Y\in
ANE(g,A,X)$ if and only if $g$ extends to a map
$\overline{g}\colon X\to Cone(Y)$.

Everywhere below $L$ always denotes a $CW$-complex.

We say that $L$ is an absolute extensor of $X$, notation $L\in
AE(X)$, if $L\in AE(g,A,X)$ for every closed $A\subset X$ and
every map $g\colon A\to L$ with $L\in ANE(g,A,X)$. We say also
that the extension dimension of $X$ is not greater than $L$,
notation $\e X\leq L$, if $L\in AE(X)$. Using Dydak's version of
the Homotopy Extension Theorem \cite[Theorem 13.7]{dy1} one can
show that if $L_1$ is homotopy equivalent to $L_2$, then   $\e
X\leq L_1$ is equivalent to $\e X\leq L_2$ for any space $X$.
Moreover, our definition of $\e$ coincides with that one of
Chigogidze \cite{ch} in case $L$ is countable and with the original definition
of Dranishnikov \cite{dr1} when compact spaces are considered.

A pair of spaces $K\subset P$ is called $L$-connected if whenever
$A\subset X$ is a closed subset of a space $X$  with $\e X\leq L$,
then  every map $g\colon A\to K$ has an extension
$\overline{g}\colon X\to P$ provided  $A$ is normally placed in
$X$ with respect to $(g,P)$.  The notion of a normally placed set
was introduced in  \cite{ch} under different notation and means
that for every continuous function  $h$ on $P$ the function
$h\circ g$ can be continuously extended over $X$. Obviously, this
condition is satisfied for every normal space $X$ and every map
$g\colon A\to K$ with $A\subset X$ closed.  We sometimes say that
a pair $K\subset P$ is $L$-connected with respect to a given class
of spaces $\mathcal B$ if the additional requirement
$X\in{\mathcal B}$ is imposed in the above definition.

Quasi-finite $CW$-complexes  were introduced in \cite{k:03}  as
countable complexes $L$ satisfying the following condition: every
finite subcomplex $K$ of $L$ is contained in a finite subcomplex
$P\subset L$  such that the pair $K\subset P$ is $L$-connected
with respect to Polish spaces.  It was also shown in \cite{k:03}
that there exists a countable quasi-finite complex $M$ extension
type $[M]$ of which does  not contain a finitely dominated complex
(see \cite{ch1} for more information on extension types). In this
note we extend the above definition by considering not necessarily
countable complexes. Here is our revised definition: a
$CW$-complex $L$ is quasi-finite if every finite subcomplex $K$ of
$L$ is contained in a finite subcomplex $P\subset L$  such that
the pair $K\subset P$ is $L$-connected. It is easy to verify that
this definition coincides with the definition given in \cite{k:03}
in case $L$ is countable.

We say that a map $f\colon X\to Y$ is $L$-invertible if for any
map $g\colon Z\to Y$  with $\e Z\leq L$ there is a map $h\colon
Z\to X$ such that $g=f\circ h$. If, in addition, $Z$ is required
to be from a given class of spaces $\mathcal B$, then we say that
the map $f$ is $L$-invertible with respect to the class $\mathcal
B$. Everywhere below $w(X)$ denotes the weight of the space $X$
and $\uin^{\tau}$ denotes Tychonov cube of weight $\tau$.

\begin{thm} The following conditions are equivalent for any $CW$-complex $L$ and an infinite cardinal $\tau$:
\begin{itemize}

\item[(1)] $L$ is quasi-finite.

\item[(2)] $\e \beta X\leq L$ whenever $X$ is a space with $\e X\leq L$.

\item[(3)] There exists an $L$-invertible map $f\colon
Y_{\tau}\to\uin^{\tau}$ such
that $Y_{\tau}$ is a compact space of
weight $\leq\tau$ and $\e Y_{\tau}\leq L$.

\item[(4)] For every $L$-connected pair $K\subset M$, where $K$ is
a compactum
of weight $\leq\tau$ and $M$ an arbitrary space,
there exists a compactum $P\subset M$ containing $K$ such that
$w(P)\leq\tau$ and the pair   $K\subset P$ is $L$-connected.
\end{itemize}
\end{thm}

\begin{proof} $(1)\Rightarrow (2)$ Suppose $\e X\leq L$ and let $f\colon A\to L$, where
$A$ is a closed subset of $\beta X$. It is well known that every
$CW$-complex is an absolute neighborhood extensor for the class of
compact spaces, so $L\in ANE(f,A,\beta X)$ and there exists a
closed neighborhood $B$ of $A$ in $\beta X$ and a map $g\colon
B\to L$ extending $f$.   Because $g(B)$ is compact, it is
contained in a finite subcomplex $K$ of $L$.  Since $L$ is
quasi-finite,  there exists a finite subcomplex $P$ of $L$ such
that the pair $K\subset P$ is $L$-connected.  We can assume that
$B$ is a zero-set in $\beta X$. Then $B\cap X$, being a non-empty
zero-set in $X$, is  normally placed in $X$ with respect to
$(g,P)$. Therefore,  the map $g\colon B\cap X\to K$ extends to a
map $h\colon X\to P$ because $\e X\leq L$  and the pair $K\subset
P$ is $L$-connected. Finally, let $\overline{h}\colon\beta X\to P$
be the unique extension of $h$. Then $\overline{h}$ extends $f$,
so $\e\beta X\leq L$.

$(2)\Rightarrow (3)$ We consider the family of all
maps
$\{h_{\alpha}\colon X_{\alpha}\to\uin^{\tau}\}_{
\alpha\in\Lambda}$ such that each $X_{\alpha}$ is a closed subset
of $\uin^{\tau}$  with $\e X_{\alpha}\leq L$.  Let $X$ be the
disjoint sum of all  $X_{\alpha}$ and the map $h\colon X\to
\uin^{\tau}$ coincides with  $h_{\alpha}$ on every $X_{\alpha}$.
Clearly, $\e X\leq L$. Therefore, $\e\beta X\leq L$. Consider the
extension $\overline{h}\colon\beta X\to\uin^{\tau}$.  Then, by the
factorization theorem from \cite{lrs}, there exists a compact
space $Y_{\tau}$ of weight $\leq\tau$ and maps $r\colon\beta X\to
Y_{\tau}$ and $f\colon Y_{\tau}\to\uin^{\tau}$ such that $\e
Y_{\tau}\leq L$ and $f\circ r=\overline{h}$.

Let us show that $f$ is $L$-invertible. Take a space $Z$ with $\e
Z\leq L$ and a map $g\colon Z\to\uin^{\tau}$. Considering  $\beta
Z$  and the  extension $\overline{g}\colon\beta Z\to\uin^{\tau}$
of $g$, we can assume that  $Z$ is always compact.  We also can
assume that the weight of $Z$ is $\leq\tau$ ( otherwise we apply
again the factorization theorem from \cite{lrs} to find a compact
space $T$ of weight $\leq\tau$ and maps $g_1\colon Z\to T$ and
$g_2\colon T\to\uin^{\tau}$ with $\e T\leq L$ and $g_2\circ
g_1=g$, and  then consider the space $T$ and the map $g_2$
instead, respectively, of $Z$ and $g$).  Therefore, without losing
generality, we can assume that  $Z$ is a closed subset of
$\uin^{\tau}$. According to the definition of $X$ and the map $h$,
there is an index $\alpha\in\Lambda$ such that $Z=X_{\alpha}$ and
$g=h_{\alpha}$.  The restriction $r|Z\colon Z\to Y_{\tau}$ is a
lifting of $g$, i.e. $f\circ (r|Z)=g$.

$(3)\Rightarrow (4)$ Suppose that $K$ is a compact  subset of the
space $M$ with $w(K)\leq\tau$ and $K\subset M$ being
$L$-connected. We embed $K$ in $\uin^{\tau}$  and consider  an
$L$-invertible mapping $f\colon Y_{\tau}\to\uin^{\tau}$  such that
$Y_{\tau}$ is compact and $\e Y_{\tau}\leq L$.  Let
$\tilde{K}=f^{-1}(K)$ and $h=f|\tilde{K}$. Obviously, $\tilde{K}$
is normally placed in $Y_{\tau}$ with respect to $(h,M)$.
Consequently, $h$ extends to a map $\overline{h}\colon Y_{\tau}\to
M$ and let $P=\overline{h}(Y_{\tau})$.  Obviously, $w(P)\leq\tau$,
so that it remains only to show that $K\subset P$ is
$L$-connected. For this end,  let $g\colon A\to K$, where
$A\subset X$ is a closed normally placed subset of $X$ with
respect to $(g,P)$ and $\e X\leq L$. This implies that $A$ is
normally placed in $X$ with respect to $(g,\uin^{\tau})$. Since
$\uin^{\tau}$ is an absolute extensor, there exists an extension
$g_1\colon X\to\uin^{\tau}$ of $g$.  Next, we lift $g_1$ to a map
$g_2\colon X\to Y_{\tau}$  such that $f\circ g_2=g_1$ (recall that
$f$ is $L$-invertible) and  let $\overline{g}=\overline{h}\circ
g_2$.  Clearly, $\overline{g}$ is a map from $X$ into $P$
extending $g$. Therefore, $K\subset P$ is $L$-connected.

$(4)\Rightarrow (1)$ Take a finite subcomplex $K$  of $L$. Let us
first show that the pair $K\subset L$ is $L$-connected. Suppose
$Z$ is a space with $\e Z\leq L$, $A\subset Z$ closed and $g\colon
A\to K$ a map such that $A$ is normally placed in $Z$ with respect
to $(g,L)$.  Since $K$  is $C$-embedded in $L$, $A$ is normally
placed in $Z$ with respect to $(g, K)$.  The last condition
together with the fact that $K$ is an absolute neighborhood
extensor for all separable metric spaces implies that $K\in ANE(g,
A,Z)$. Indeed, we embed  $K$ in $\R^{\omega}$  and fix a
retraction $r\colon U\to K$, where $U$ is a neighborhood of $K$ in
$\R^{\omega}$ . Since  $A$ is normally placed in $Z$ with respect
to $(g, K)$, we can find a map $h\colon Z\to \R^{\omega}$
extending $g$. Then $h^{-1}(U)$ is a co-zero neighborhood of $A$
in $Z$ which contains the zero-set $h^{-1}(K)$  and $r\circ
h\colon h^{-1}(U)\to K$  extends $g$. Hence,  $K\in ANE(g, A,Z)$
which yields $L\in ANE(g, A,Z)$.  Since  $\e Z\leq L$, $g$ can be
extended to a map $\overline{g}\colon Z\to L$. Thus, $K\subset L$
is an $L$-connected pair. Therefore there exists a compact set
$H\subset L$ containing $K$ such that the pair $K\subset H$ is
$L$-connected. Finally, we take a finite subcomplex $P$ of $L$
which contains $H$ and observe that the pair $K\subset P$ is also
$L$-connected.  Hence, $L$ is quasi-finite.
\end{proof}

\begin{cor} None of the Eilenberg-MacLane complexes $K(G,n)$, $n\geq 2$ and $G$ an Abelian group, is quasi-finite.
\end{cor}

\begin{proof}
This follows from Theorem 2.1(2) and the following statement (see
\cite[Theorem 1.4]{l1}): there exists a separable metric space $X$
with $\dim_G X\leq 2$ and $\e\beta X>L$ for every Abelian group
$G$ and every non-contractible $CW$-complex $L$. Here $\dim_G X$
denotes the cohomological dimension of $X$ with respect to the
group $G$.
\end{proof}

Let us also observe that for every quasi-finite complex $L$ there
exists a compact metrizable space which is universal for the class
of all separable metric spaces of $\e\leq L$, in particular every
space from this class has a compactification of $\e\leq L$.
Indeed, let $Y_{\omega}$ be the space from Theorem 2.1(3). Then,
for every $X$ from the above class we take an embedding $i\colon
X\to\uin^{\omega}$ and lift $i$ to a map $j\colon X\to
Y_{\omega}$.  The required compactification of $X$ is the closure
of $j(X)$ in $Y_{\omega}$.  Next corollary provides a
characterization of quasi-finite countable complexes  in terms of
compactifications.

\begin{cor}For a countable complex $L$ the following conditions are equivalent:

\begin{itemize}
\item[(a)] $L$ is quasi-finite.

\item[(b)] For every separable metrizable space $X$ with $\e X\leq
L$ and its metrizable compactification $c(X)$ there exists a
metrizable compactification $c^{*}(X)$ such that $\e c^*(X)\leq L$
and $c^{*}(X)\geq c(X)$ $($i.e., there is a map from $c^{*}(X)$
onto $c(X)$ which is the identity on $X$$)$.
\end{itemize}
\end{cor}

\begin{proof}$(a)\Rightarrow (b)$
Let $L$ be quasi-finite and $X$ a separable metric space with $\e
X\leq L$.  We take a metric compactification $c(X$) of $X$ and a
map $f\colon \beta X\to c(X)$ such that $f(x)=x$ for every $x\in
X$. Since, by Theorem 2.1, $\e\beta X\leq L$, $f$ can be factored
through a metrizable compactum $Z$ with $\e Z\leq L$. Clearly, $Z$
is a compactification of $X$ which is $\geq c(X)$.

$(b)\Rightarrow (a)$ According to \cite[Corollary 3.4]{dm}, there
exists a metrizable compactum $Y$ with $\e Y\leq L$ and a
surjective map $f\colon Y\to\uin^{\omega}$  such that for any map
$g\colon X\to\uin^{\omega}$,   $X$ being separable metrizable with
$\e X\leq L$, there exists an embedding $i\colon X\to Y$ lifting
$g$, i.e.  $f\circ i=g$.  Hence, $f$ is $L$-invertible with
respect to separable metric spaces.  By Theorem 2.1(3), it
suffices to show that $f$ is $L$-invertible.  Consider $g\colon
Z\to\uin^{\omega}$ where $\e Z\leq L$.  According to
\cite[Proposition 4.9]{ch}, there exist a Polish space $P$ with
$\e P\leq L$ and maps $h\colon Z\to P$ and $q\colon
P\to\uin^{\omega}$ with $g=q\circ h$. We lift $q$ to a map
$\overline{q}\colon P\to Y$ such that $f\circ\overline{q}=q$. Then
$\overline{q}\circ h$ is the required lifting of $g$.
\end{proof}

Here is another property of quasi-finite complexes:

\begin{pro}   Every quasi-finite complex $L$ has the following connected-pairs property:

\begin{itemize}
\item [$(CP)$] For any metrizable compactum $K$ with $\e K\leq L$
there exists a metrizable compactum $P$ containing $K$ such that
$\e P\leq L$ and the pair $K\subset P$ is $L$-connected.
\end{itemize}
\end{pro}

\begin{proof} Suppose $K$
is a metrizable compactum with $\e K\leq L$. We embed $K$ into the
Hilbert cube $\uin^{\omega}$ and take an $L$-invertible map
$f\colon Y\to\uin^{\omega}$ such that $Y$ is a metrizable
compactum with $\e Y\leq L$ (see Theorem 2.1(3)).  Consider the
adjunction space $Y\cup_fK$, i.e. the disjoint union of
$Y-f^{-1}(K)$ and $K$ with the topology consisting of the usual
open subsets of  $Y-f^{-1}(K)$ together with sets of the form
$f^{-1}(U-K)\cup(U\cap K)$ for open subsets $U$ of
$\uin^{\omega}$. There are two associated maps $p_K\colon Y\to
Y\cup_fK$ and $f_K\colon Y\cup_fK\to \uin^{\omega}$ such that
$f=f_K\circ p_K$.  Since $f$ is $L$-invertible, so is $f_K$.
Moreover, $Y-f^{-1}(K)$, being open in $Y$, is the  union of
countably many compact sets each with $\e\leq L$.  Hence, by the
countable sum theorem, $\e Y\cup_fK\leq L$.

We need only to show that the pair $K\subset Y\cup_fK$ is
$L$-connected. Let $g\colon A\to K$  be a map from a closed subset
$A\subset Z$ such that $\e Z\leq L$ and  $A$ is normally placed in
$Z$ with respect to $(g,Y\cup_fK)$. Then, considering $g$ as a map
from $A$ into $K\subset\uin^{\omega}$,   we obviously have that
$A$ is normally placed in $Z$ with respect to $(g,
\uin^{\omega})$.  Since $ \uin^{\omega}$ is an absolute extensor,
there exists a map $\overline{g}\colon Z\to \uin^{\omega}$
extending $g$. Finally, since $f_K$ is $L$-invertible, we lift
$\overline{g}$ to a map  $h\colon Z\to Y\cup_fK$ with $f_K\circ
h=\overline{g}$. Clearly, $h$ extends $g$.
\end{proof}

\begin{pro}
For every $n\geq 2$ there is no $K(\mathbb{Z},n)$-connected pair
$K\subset P$ of compact sets such that $K$ is homeomorphic to the
$n$-dimensional sphere  ${\bbb S}^n$ and $\dim_ {\mathbb{Z}}P\leq
n$.
\end{pro}

\begin{proof}
We use the arguments from the proof of \cite[Theorem 3.5]{dm}.
Suppose for some $n\geq 2$ there is a $K(\mathbb{Z},n)$-connected
compact pair ${\bbb S}^n\subset P$ with $\dim_ {\mathbb{Z}}P\leq
n$.   We choose a complex $L$ of type $K(\mathbb{Z},n)$ and having
finite skeleta. It was shown in \cite{dw} that there exist
metrizable compacta $X_k$, $k\geq 1$, such that:

\begin{itemize}
\item $\dim_ {\mathbb{Z}}X_k\leq n$ for each $k$;
\item each $X_k$
contains a copy of ${\bbb S}^n$; \item the inclusion $i\colon{\bbb
S}^n \hookrightarrow L$ cannot be extended over $X_k$ so that the
image of the extension is contained in the $k$-skeleton $L^{(k)}$
of $L$.
\end{itemize}

We take an extension $h\colon P\to L$ of the inclusion
$i\colon{\bbb S}^n \hookrightarrow L$, and $m$ such that $h(P)\subset L^{(m)}$.  This means that the inclusion
$j\colon{\bbb S}^n \hookrightarrow P$ cannot be extended to a map from $X_m$ into $P$ which contradicts the fact that
${\bbb S}^n\subset P$ is $L$-connected.
\end{proof}

The problem \cite{we} whether, for any fixed $n\geq 2$  there is a
universal space in the class of all metrizable compacta $X$  with
$\dim_{\mathbb{Z}}\leq n$ is still unsolved. Zarichnyi \cite{z}
observed that each of the above classes does not have an universal
element  which is an absolute extensor for the same class.
Proposition 2.5 yields a little bit stronger observation.

\begin{cor} None of the  complexes $K(\mathbb{Z},n)$, $n\geq 2$, have the $(CP)$-property.
\end{cor}
Recall that a map $f\colon X\to Y$ between metrizable spaces is
called uniformly 0-dimensional \cite{mk} if there exists a metric
on $X$ generating its topology such that for every $\epsilon>0$
every point of $f(X)$ has a neighborhood $U$ in $Y$ with
$f^{-1}(U)$ being the union of disjoint open subsets of $X$ each
of diameter $<\epsilon$.  It is well known that every metric space
admits uniformly 0-dimensional map into $l_2$.

\begin{pro} Let $L$ be a quasi-finite $CW$-complex.
Then for every $\tau\geq\omega$ there exists  a perfect
$L$-invertible surjection $f_{(L,\tau)}\colon Y_{(L,\tau)}\to
l_2(\tau)$ such that:

\begin{itemize}
\item[(a)] $Y_{(L,\tau)}$ is a completely metrizable space of
weight $\tau$ with $\e Y_{(L,\tau)}\leq L$.

\item[(b)] Every $($completely$)$ metrizable space of weight
$\leq\tau$ and extension dimension $\leq L$ can be embedded as a
$($closed$)$ subspace of $Y_{(L,\tau)}$.
\end{itemize}
\end{pro}

\begin{proof}
By Theorem 2.1(3), there exists an $L$-invertible map $f\colon
Y\to\uin^{\omega}$, where $Y$ is a metrizable compactum with $\e
Y\leq L$.  We embed $l_2$ in $\uin^{\omega}$ and let
$Y_{(L,\omega)}=f^{-1}(l_2)$ and
$f_{(L,\omega)}=f|Y_{(L,\omega)}$. Then $\e Y_{(L,\omega)}\leq L$
and since $f$ is $L$-invertible, so is $f_{(L,\omega)}$.

If $\tau >\omega$, we take a metric $d_1$ on $l_2(\tau)$ and a
uniformly 0-dimensional map $g\colon l_2(\tau)\to l_2$ with
respect to $d_1$. Denote by $Y_{(L,\tau)}$  the fibered product of
$l_2(\tau)$ and $Y_{(L,\omega)}$ with respect to the maps $g$ and
$f_{(L,\omega)}$. We also consider the projections
$f_{(L,\tau)}\colon Y_{(L,\tau)}\to l_2(\tau)$ and  $h\colon
Y_{(L,\tau)}\to Y_{(L,\omega)}$. Since $f_{(L,\omega)}$ is a
perfect and $L$-invertible surjection, so is $f_{(L,\tau)}$.  If
$d_2$  is any metric on $Y_{(L,\omega)}$, then $h$ is uniformly
0-dimensional with respect to the metric $\displaystyle
d=\sqrt{d_1^2+d_2^2}$ on $Y_{(L,\tau)}$ (see \cite{ap}). Thus
$Y_{(L,\tau)}$ admits a uniformly 0-dimensional map into the space
$Y_{(L,\omega)}$ having extension dimension $\leq L$.  Hence, by
\cite[Theorem 1.2]{ml}, $\e Y_{(L,\tau)}\leq L$.  Observe that
$Y_{(L,\tau)}$ is completely metrizable as a perfect preimage of
the completely metrizable space $l_2(\tau)$.

To prove the second item,  suppose $M$ is a metrizable space of
weight $\leq\tau$ and $\e M\leq L$. We consider $M$ as a subset of
$l_2(\tau)$ and use the $L$-invertibility of $f_{(L,\tau)}$ to
lift the identity map on $M$. Obviously this lifting is an
embedding  of $M$ into $Y_{(L,\tau)}$. Moreover, if $M$ is
completely metrizable, then we can embed it in $l_2(\tau)$ as a
closed subspace. This implies that the corresponding embedding of
$M$ in $Y_{(L,\tau)}$ is also closed.
\end{proof}

A completion theorem for $L$-dimensional metric spaces, where $L$
is any countable $CW$-complex, was established in \cite{wo}.  It
follows from Proposition 2.7 that this is also true for
quasi-finite (not necessarily countable) complexes $L$.

\begin{cor} Let $L$ be a quasi-finite complex.
Then every metrizable space $X$ with $\e X\leq L$ has a completion
with extension dimension $\leq L$.
\end{cor}

\begin{cor} Let $L$ be a quasi-finite complex and $X$
a metrizable space. Then $\e X\leq L$ if and only if $X$ admits a
uniformly 0-dimensional map into a separable metrizable space of
extension dimension  $\leq L$.
\end{cor}

\begin{proof} In one direction (sufficiency)
this follows from the mentioned above result of Levin
\cite[Theorem 1.2]{ml}.  Suppose $X$ is a metrizable space of
weight $\tau$ with $\e X\leq L$. By    Proposition 2.7, $X$ can be
embedded in the space $Y_{(L,\tau)}$. It follows from the
construction of $Y_{(L,\tau)}$ that the map $h\colon
Y_{(L,\tau)}\to Y_{(L,\omega)}$ is uniformly 0-dimensional. Then
the restriction $h|X$ is also uniformly 0-dimensional which
completes the proof.
\end{proof}

A general factorization theorem for $L$-dimensional compact
spaces, where $L$ is an arbitrary complex, was proved in
\cite{lrs}. We provide here a factorization theorem for
$L$-dimensional metrizable spaces with $L$ being quasi-finite (see
\cite[Theorem 1.5]{ml} for similar result with  $L$ countable).

\begin{pro} Let $L$ be a quasi-finite complex and let $f\colon X\to Y$ be a map with $Y$ metrizable.
If $\e X\leq L$, then $f$ factors through a metrizable space $Z$ such that $\e Z\leq L$ and  $w(Z)\leq w(Y)$.
\end{pro}

\begin{proof} Let us first show how to reduce this proposition to the case $Y$ is separable.
This reduction is well known (see, for example, \cite{ap}), but we
present it here for the reader's convenience.   Suppose the result
holds when the range space is separable and metrizable. We take a
uniformly 0-dimensional map $g\colon Y\to l_2$ and apply the
``separable factorization theorem" to the map $g\circ f\colon X\to
l_2$ to obtain a separable metrizable space $M$ and maps $q\colon
X\to M$ and $h\colon M\to l_2$ with $\e M\leq L$ and $h\circ
q=g\circ f$. Let $p_M\colon Z\to M$ and $p_Y\colon Z\to Y$ be the
pullbacks of $g$ and $h$ respectively. Clearly, $Z$ is a
metrizable space of weight $w(Z)\leq w(Y)$. Since $g$ is uniformly
0-dimensional, so is $p_M$. Then, by \cite[Theorem 1.2]{ml}, $\e
Z\leq L$.

Now we prove the ``separable case".  Let $\tilde{Y}$ be a
metrizable compactification of $Y$ and $\tilde{f}\colon \beta
X\to\tilde{Y}$ be the \v{C}ech-Stone extension of $f$.  Since $L$
is quasi-finite, $\e\beta X\leq L$. Therefore we can apply the
factorization theorem of Levin-Rubin-Schapiro \cite{lrs} to obtain
a metrizable compactum $\tilde{Z}$ and maps
$\tilde{f_1}\colon\beta X\to\tilde{Z}$ and
$\tilde{f_2}\colon\tilde{Z}\to\tilde{Y}$ such that
$\tilde{f_2}\circ\tilde{f_1}=\tilde{f}$ and $\e\tilde{Z}\leq L$.
Then the space $Z=\tilde{f_1}(X)$
 and the maps $f_1=\tilde{f_1}|X$ and $f_2=\tilde{f_2}|Z$ form the required factorization.
 \end{proof}

We say that a map $f\colon X\to Y$ is $L$-soft, where $L$ is a
$CW$-complex, if for any space $Z$ with $\e Z\leq L$, any closed
set $A\subset Z$ and any two maps $h\colon Z\to Y$ and $g\colon
A\to X$, where $A$ is normally placed in $Z$ with respect to
$(g,X)$ and $f\circ g=h|A$, there exists a map $\overline{g}\colon
Z\to X$ extending $g$ such that $f\circ\overline{g}=h$.  If, in
the above definition, we additionally require $Z$ to be from a
given class of spaces $\mathcal A$, then we say  that $f$ is
$L$-soft with respect to the class $\mathcal A$. It was
established in \cite{chv} that for every countable complex $L$ and
every metric space $Y$ there exists an $L$-soft map $f\colon X\to
Y$ such that $X$ is a metric space of extension dimension $\leq L$
and $w(X)=w(Y)$. We are going to show that quasi-finite complexes
also have this property.

\begin{pro} Let $L$ be a quasi-finite $CW$-complex.
Then for every $\tau\geq\omega$ there exists  an $L$-soft map $p_{(L,\tau)}\colon X_{(L,\tau)}\to l_2(\tau)$ such that:

\begin{itemize}
\item[(a)] $X_{(L,\tau)}$ is a completely metrizable space of weight $\tau$ with $\e X_{(L,\tau)}\leq L$.

\item[(b)] $X_{(L,\tau)}$ is an absolute extensor for all metrizable spaces of $\e \leq L$.

\item[(c)] $p_{(L,\tau)}$ is a strongly $(L,\tau)$-universal map,
i.e. for any open cover $\mathcal U$ of $X_{(L,\tau)}$, any
$($complete$)$ metrizable space $Z$ of weight $\leq\tau$ with $\e
Z\leq L$ and any map $g\colon Z\to X_{(L,\tau)}$ there exists a
$($closed$)$ embedding $h\colon Z\to X_{(L,\tau)}$ which is
$\mathcal U$-close to $g$ and $p_{(L,\tau)}\circ
g=p_{(L,\tau)}\circ h$.
\end{itemize}
\end{pro}

\begin{proof}
Using Proposition 2.11 and following Zarichnyi's idea from
\cite{z} (see also \cite{ch}) that invertibility generates
softness, we can show the existence of a complete separable
metrizable space $X$ with $\e X\le L$ and an $L$-soft map $f\colon
X\to l_2$. Then, as in \cite{chv}, we construct the space
$X_{(L,\tau)}$ and the map $p_{(L,\tau)}\colon X_{(L,\tau)}\to
l_2(\tau)$ possessing the desired properties.
\end{proof}


\section{Some more properties of quasi-finite complexes}

In this section, all spaces and all $CW$-complexes, unless stated
otherwise, are, respectively, metrizable and quasi-finite. We are
going to show that some properties of finitely dominated complexes
remain valid for quasi-finite complexes. We say that a space $X$
is an absolute (neighborhood) extensor in dimension $L$ (notation
$X\in A(N)E(L)$) if for every space $Z$ of extension dimension
$\leq L$ and every map $g\colon A\to X$, where $A$ is a closed
subset of $Z$,  there exists an extension of $g$ over $Z$ (resp.,
over a neighborhood of $A$ in $Z$).

Everywhere below $cov(X)$ denotes the family of all open covers of
$X$. Two maps $f_0, f_1\colon X\to Y$ are $L$-homotopic \cite{ch1}
if for any map $h\colon Z\to X\times [0,1]$, where $Z$ is a  space
with $\e Z\leq L$, the composition $(f_0\oplus f_1)\circ
h|\big(h^{-1}\big(X\times\{0,1\}\big)\big)\colon
h^{-1}\big(X\times\{0,1\}\big)\to Y$ admits an extension $H\colon
Z\to Y$. If $\mathcal U\in cov(X)$  and the extension $H$ in the
above definition can be chosen such that the collection
$\{H\big(h^{-1}(\{x\}\times [0,1])\big): x\in X\}$ refines
$\mathcal U$, then $f_0$ and $f_1$ are called $({\mathcal
U},L)$-homotopic.

The following three propositions were given in \cite{ch1} for
finitely dominated countable complexes $L$ and Polish
$ANE(L)$-spaces $X$. Because of Proposition 2.7, one can show they
also hold for quasi-finite complexes $L$ and arbitrary (not
necessarily Polish) $ANE(L)$-spaces.

\begin{pro}
Let $X$ be an $ANE(L)$-space and ${\mathcal U}\in cov(X)$. Then
there exists a cover ${\mathcal V}\in cov(X)$ such that any two
$\mathcal V$-close maps of any space into $X$ are $({\mathcal
U},L)$-homotopic.
\end{pro}

\begin{pro}
Let  $X\in ANE(L)$ and ${\mathcal U}\in cov(X)$. Then there exists
a cover ${\mathcal V}\in cov(X)$ refining $\mathcal U$, such that
the following condition holds:

\begin{itemize}
\item[(H)] For any space $Z$ with $\e Z\leq L$, any closed
$A\subset Z$, and any two $\mathcal V$-close maps $f,g\colon A\to
X$ such that $f$ has an extension $F\colon Z\to X$, it follows
that $g$ also can be extended to a map $G\colon Z\to X$ which is
$({\mathcal U},L)$-homotopic to $F$.
\end{itemize}
\end{pro}

\begin{pro}
Let $X\in ANE(L)$,  $Z$ be a space with $\e Z\leq L$ and $A\subset
Z$ closed. If $f,g\colon A\to X$ are $L$-homotopic and $f$ admits
an extension $F\colon Z\to X$, then  $g$ also admits an extension
$G\colon Z\to X$, and we may be assume that $F$ and $G$ are
$L$-homotopic.
\end{pro}

A pair of closed subsets $X_0\subset X_1$ of a space $X$ is called
$UV(L)$-connected in $X$ if every neighborhood $U$ of $X_1$ in $X$
contains a neighborhood $V$ of $X_0$ such that $V\subset U$ is
$L$-connected with respect to metrizable spaces, i.e. any map
$g\colon A\to V$, where $A$ is a closed subset of a space $Z$ with
$\e \leq L$, admits an extension $\overline{g}\colon Z\to U$. When
$X_0\subset X_0$ is $UV(L)$-connected in $X$, we say that $X_0$ is
$UV(L)$ in $X$.  If in the above definition  all pairs under
consideration are $L$-connected with respect to a given class
$\mathcal A$, we obtain the notion of $UV(L)$-sets with respect to
$\mathcal A$. If instead of $L$-connectedness of the pair
$V\subset U$ we require the inclusion $V\subset U$ to be
$L$-homotopic to a constant map in $U$ then the pair $X_0\subset
X_1$ (resp. the set $X_0$) is called $UV(L)$-homotopic in $X$.
Obviously, every $UV(L)$-connected pair is $UV(L)$-homotopic. Next
corollary, which follows from Proposition 3.3,  shows that these
two properties are equivalent in case $X\in ANE(L)$.

\begin{cor}
Let $X$ be an $ANE(L)$-space. A pair $X_0\subset X_1$ of closed
subsets of $X$ is $UV(L)$-connected in $X$ if and only if it is
$UV(L)$-homotopic in $X$.
\end{cor}

\begin{lem}
Let $X_0\subset X_1\subset X\subset E$, where both $X$ and $E$ are
$ANE(L)$-spaces and $X\subset E$  is closed. Then  the pair
$X_0\subset X_1$ is $UV(L)$-connected in $X$ if and only if it is
$UV(L)$-connected in $E$.
\end{lem}

\begin{proof} By Proposition 2.7, there exists a perfect
$L$-invertible surjection $f\colon\tilde{E}\to E$ with $\e
\tilde{E}\leq L$, and let $\tilde{X}=f^{-1}(X)$. Since $X\in
ANE(L)$, we can extend $f|\tilde{X}$ to a map $g\colon W\to X$
with $W$ being a neighborhood of $\tilde{X}$ in $\tilde{E}$. Since
$f$ is closed, we may assume that $W=f^{-1}(G)$ for some
neighborhood $G$ of $X$ in $E$. The claim below follows from our
constructions.

{\em Claim.} For every open $O\subset X$ the set
$O^*=G-f\big(g^{-1}(X-O)\big)$ is open in $G$ and has the
following two properties: $O^*\cap X=O$ and
$g\big(f^{-1}(O^*)\big)=O$.

Suppose $X_0\subset X_1$ is $UV(L)$-connected in $X$. We are going
to show that this pair is $UV(L)$-connected in $E$. To this end,
let $U\subset G$ be a neighborhood of $X_1$ in $E$. Then there is
a neighborhood $O$ of $X_0$ in $X$ such that $O\subset U\cap X$ is
$L$-connected. Since $U$ is an $ANE(L)$ (as an open subset of
$E$), we can apply Proposition 3.2 for the space $U$ and the
one-element cover $\mathcal U=\{U\}$ to find an open cover
$\mathcal V=\{V_{\alpha}:\alpha\in\Lambda\}$ of $U$ satisfying the
condition $(H)$. For every $\alpha$ let
$G_{\alpha}=V_{\alpha}\cap\big(V_{\alpha}\cap X\big)^*\cap O^*$
and $V=\bigcup\{G_{\alpha}: \alpha\in\Lambda\}$. Obviously,
$V\subset U$ is open and contains $X_0$. The pair $V\subset U$ is
$L$-connected. Indeed, let $h\colon A\to V$ be a map, where
$A\subset Z$ is closed and $\e Z\leq L$. Since $f$ is
$L$-invertible, $h$ admits a lifting $h_1\colon A\to f^{-1}(V)$,
i.e. $h=f\circ h_1$. According to the Claim,
$g\big(f^{-1}(G_{\alpha})\big)\subset V_{\alpha}\cap X$,
$\alpha\in\Lambda$, and $V\cap X\subset O$. This implies that $h$
and the map $h_2=g\circ h_1\colon A\to V\cap X$ are $\mathcal
V$-close. Since the pair $O\subset U\cap X$ is $L$-connected,
$h_2$ can be extended to a map from $Z$ into $U\cap X$.   This
yields, according to Proposition 3.2,   that $h$ also can be
extended to a map from $Z$ into $U$.

Now,  suppose the pair $X_0\subset X_1$ is $UV(L)$-connected in
$E$.  To show this pair is $UV(L)$-connected in $X$,  let  $U$ be
a neighborhood of $X_1$ in $X$. Then $U^*\subset G$ is open in
$E$, and we can find a neighborhood $V$ of $X_0$ in $E$ such that
$V\subset U^*$ is $L$-connected.  The pair $V\cap X\subset U$ is
$L$-connected. Indeed, any map  $h\colon A\to V\cap X$, where
$A\subset Z$ is closed and $\e Z\leq L$, admits an extension
$h_1\colon Z\to U^*$. Then the map $\overline{h}=g\circ h_2\colon
Z\to U$, where $h_2\colon Z\to f^{-1}(U^*)$ is a lifting of $h_1$,
extends $h$.
\end{proof}

\begin{thm} Suppose  $X$ is an $ANE(L)$-space and the pair
$X_0\subset X_1$ is $UV(L)$-connected in  $X$. Then it is
$UV(L)$-connected in any $ANE(L)$-space in which $X_1$ is
embeddable as a closed subspace.
\end{thm}

\begin{proof}
Let   $i\colon X_1\to Y$ be a closed embedding, where $Y\in
ANE(L)$,  and $M$ be the space obtained from the disjoint union
$X\uplus Y$ by identifying all pairs of points $x\in X_1\subset X$
and $i(x)\in Y$.  The space $M$ is metrizable and if $p\colon
X\uplus Y\to M$ is the quotient map, then $p(X)$, $p(Y)$ and
$p(X_1)$ are closed sets in $M$ homeomorphic, respectively, to
$X$,  $Y$ and $X_1$. Moreover, $p(X_1)$ is the common part of
$p(X)$ and $p(Y)$. We embed $M$  in a normed space $E$ as a closed
subspace. Every normed space is an absolute extensor for the class
of metrizable spaces, so $E\in ANE(L)$. Since the pair
$p(X_0)\subset p(X_1)$ is $UV(L)$-connected in $p(X)$, by Lemma
3.5 it is also $UV(L)$-connected in $E$.  This implies, again by
Lemma 3.5,  that $p(X_0)\subset p(X_1)$ is $UV(L)$-connected in
$p(Y)$.
\end{proof}

\begin{cor}
If a space $X$ is $UV(L)$ in a given $ANE(L)$-space,  then $X$ is
$UV(L)$ in any $ANE(L)$-space in which $X$ is embeddable as a
closed subset.
\end{cor}

In the existing literature, the $UV^n$- property, and more
general, the $UV(L)$-property, is defined for compact spaces, see
\cite{ch1} and \cite{bchk}. We extend this definition to arbitrary
(metrizable) spaces: $X$ is a {\em $UV(L)$-space} if it is $UV(L)$
in some $ANE(L)$-space containing $X$ as a closed subspace.
According to Corollary 3.7, the $UV(L)$-property does not depend
on the embeddings in $ANE(L)$-spaces (for compact spaces and
finite complexes $L$ this was done in \cite{bchk}). It follows
from Corollary 3.4 that $X$ is a $UV(L)$-space if and only if $X$
is $UV(L)$-homotopic in every space $Y\in ANE(L)$ containing $X$
as a closed subset.

Recall that  a normal space $X$ is a $C$-space \cite{ag} if for
any sequence $\{\omega_n\}$ of open covers of $X$ there exists a
sequence $\{\gamma_n\}$ of open disjoint families such that each
$\gamma_n$ refines $\omega_n$ and $\cup\gamma_n$ covers $X$. Every
finite-dimensional paracompactum, as well as every
countable-dimensional metrizable space has property $C$
\cite{re:95}.

We say that a complex $L$ (not necessarily quasi-finite) possesses
{\em the soft map property} if for every space $X$ there exists a
space $Y$ with $\e Y\leq L$ and an $L$-soft  map from $Y$ onto
$X$. Every countable complex has the soft map property (see
\cite{chv}), as well as every quasi-finite complex (by Proposition
2.11).

A pair of spaces $\tilde{V}\subset\tilde{U}$ is called an
$L$-extension of the pair $V\subset U$ \cite{bch} if $\tilde{U}\in
AE(L)$ and there exists a map $q\colon\tilde{U}\to U$ such that
the restriction $q|\tilde{V}$ is an $L$-soft map onto $V$.  The
following property of $L$-extension pairs was established in
\cite{bch}.

\begin{lem}
Let  $L$ be a complex (not necessarily quasi-finite) with the soft
map property and $\tilde{V}\subset\tilde{U}$  an $L$-extension of
the pair $V\subset U$. Let also  $A\subset B$ be a pair of closed
subsets of a space $X$ with $\e X\leq L$. Suppose we have maps
$f\colon B\to U$ and $g\colon A\to\tilde{U}$ such that $q\circ
g=f|A$ and $f\big(\overline{B\backslash A}\big)\subset V$. Then
there exists a map $h\colon X\to\tilde{U}$ such that $q\circ
(h|B)=f$.
\end{lem}

\begin{lem}
Let  $L$ be a complex (not necessarily quasi-finite). Every
$L$-connected pair $V\subset U$ of spaces admits an $L$-extension
provided $L$ has the soft map property.
\end{lem}

\begin{proof}
We take a normed space $E$ containing $V$ as a closed subspace and
an $L$-soft surjection $g\colon\tilde{U}\to E$ such that
$\tilde{U}$ is a space of $\e \leq L$.  Since $V\subset U$ is
$L$-connected, there exists a map $q\colon\tilde{U}\to U$
extending the map $g|\tilde{V}$, where $\tilde{V}=g^{-1}(V)$.
Moreover, $\tilde{U}\in AE(L)$ because $E$ is an absolute extensor
for the class of metrizable spaces and $g$ is $L$-soft.
Therefore, $\tilde{V}\subset\tilde{U}$  is an $L$-extension of
$V\subset U$.
\end{proof}

If $A$ is a subset of a space $X$ we denote the star of $A$ with
respect to a cover $\omega\in cov(X)$ by $\St (A,\omega )$. We say
that $\nu\in cov (X)$ is a strong star-refinement of $\omega\in
cov (X)$ if for each $V\in\nu$ there exists $W\in\omega$ such that
$\St (V,\nu)\subset W$.

{\bf Auxiliary Construction.} Suppose we are given the spaces $X$,
$Z$  and  the map  $g\colon A\to X$, where $A\subset Z$ is closed.
Let $\alpha_n=\{U_n(x):x\in X\}$, $\beta_{n}=\{V_{n}(x):x\in X\}$,
$n\geq 0$, be two sequences of open covers of $X$ and $\mu_n^*$,
$n\geq 1$, be a sequence of disjoint open families in $A$ such
that:

\begin{itemize}
\item[(1)] $\alpha_n$ is a strong star refinement of $\beta_{n-1}$
for any $n\geq 1$.

\item[(2)] each $\mu^*_n$, $n\geq 1$, refines $g^{-1}(\beta_n)$
and $\cup\{\mu_n^{*}:n\geq 1\}$ is a locally finite cover of  $A$.
\end{itemize}

We are going first to construct open and disjoint families
$\mu_n$, $n\geq 1$, in $Z$ satisfying the following condition:

\begin{itemize}
\item[(3)] $\mu=\cup\{\mu_n:n\geq 1\}$ is locally finite in  $Z$
and the restriction of each $\mu_n$ on $A$ is $\mu_n^*$.
\end{itemize}

To this end, we choose  an upper semi-continuous (br., u.s.c.)
set-valued map $r\colon Z\to A$ such that each $r(z)$ is a finite
set and $r(z)=\{z\}$ for $z\in A$ (see \cite{n} for the existence
of such $r$). Recall that $r$ is upper semi-continuous means that
$r^{\sharp}(T)=\{z\in Z: r(z)\subset T\}$ is open in $Z$ whenever
$T$ is open in $A$.    Obviously, $r^{\sharp}(T)\cap A=T$  and
$r^{\sharp}(T_1)\cap r^{\sharp}(T_2)\neq\emptyset$ if and only if
$T_1\cap T_2\neq\emptyset$ for any open subsets $T$, $T_1$ and
$T_2$ of $A$. Therefore all  families $\mu_n=\{r^{\sharp}(T):
T\in\mu_n^*\}$, $n\geq 1$,  are open and disjoint in $Z$. Since
$\mu^*$ is locally finite in $A$ and $r$ is finite-valued, the
family $\mu=\cup\{\mu_n:n\geq 1\}$ is locally finite in  $Z$.

The second part of our construction is to find points $x_W\in X$ such that

\begin{itemize}
\item[(4)] $St\big(g(W\cap A), \alpha_{n}\big)\subset V_{n-1}\big(x_W\big)$ for every $W\in\mu_n$ and $n\geq 1$
\end{itemize}

This can be done as follows.  Since $\alpha_n$ is a strong star
refinement of $\beta_{n-1}$ and $\mu_n$ refines $g^{-1}(\beta_n)$,
for every $n\geq 1$ and $W\in\mu_n$ there exist $S\in \beta_n$ and
a point $x_W\in X$ such that $St\big(g(W\cap A),
\alpha_{n}\big)\subset St\big(S,\alpha_n\big)\subset
V_{n-1}\big(x_W\big)$. The auxiliary construction is completed.

\begin{lem}
Let $L$ be a complex  $($not necessarily quasi-finite$)$ with the
soft map property and $f\colon M\to X$ be a surjection with the
following property:
\begin{itemize}
\item[(UV)]for every $x\in X$ and  its neighborhood $U(x)$ in $X$
there exists a smaller neighborhood $V(x)$ of $x$ such that the
pair $\tilde{V}(x)=f^{-1}(V(x))\subset\tilde{U}(x)=f^{-1}(U(x))$
is $L$-connected with respect to the class of metrizable spaces.
\end{itemize}
Suppose $p\colon Y\to Z$ is a surjective map with $\e Y\leq L$.
Then, for any $\omega\in cov(X)$ and any map $g\colon A\to X$,
where $A$ is a closed subset of $Z$ such that either $A$ or $g(A)$
is a $C$-space, there is a neighborhood $G$ of $A$ in $Z$ and a
map $h\colon p^{-1}(G)\to M$ with $(f\circ h)|p^{-1}(A)$ being
$\omega$-close to $g\circ p$.
\end{lem}

\begin{proof}
For every $x\in X$ and $n=0,1,2,\dots$ we choose a point $P(x)\in
f^{-1}(x)$ and neighborhoods $U_n(x)$ and $V_n(x)$ of $x$  in $X$
such that the cover $\alpha_0=\{U_0(x):x\in X\}$ refines $\omega$,
each pair $\tilde{V}_n(x)\subset\tilde{U}_n(x)$ is $L$-connected
with respect to all metrizable spaces and the covers
 $\alpha_n=\{U_n(x):x\in X\}$,  $\beta_{n}=\{V_{n}(x):x\in X\}$ satisfy condition (1) from
the auxiliary construction. Since either $A$ or $g(A)$ is a
$C$-space, there exists a sequence of disjoint open families
$\{\mu_n^{*}: n\geq 1\}$ in $A$ satisfying condition (2) above.
Therefore, according to the auxiliary construction, we can extend
each $\mu_n^*$ to a disjoint open family $\mu_n$ in $Z$ such that
$\mu=\cup\{\mu_n:n\geq 1\}$ is locally finite in  $Z$ and let $G$
be the union of all elements of $\mu$.

We introduce the following notations: $B=p^{-1}(A)$,
$\overline{g}=g\circ (p|B)$, $\Omega =p^{-1}(G)$, and
$\nu_n=p^{-1}(\mu_n)$. Obviously, each $\nu_n$ is a disjoint open
family in $Y$ and $\nu=\cup\{\nu_n:n\geq 1\}$ is a locally finite
cover of $\Omega$. Let us also consider the open covers
$\tilde{\omega}=f^{-1}(\omega)$,
$\tilde{\alpha}_n=\{\tilde{U}_n(x):x\in X\}$ and
$\tilde{\beta}_n=\{\tilde{V}_n(x):x\in X\}$ of $M$
 corresponding, respectively,  to $\omega$,
$\alpha_n$ and $\beta_n$.  According to Lemma 3.9, every pair
$\tilde{V}_{n}(x)\subset\tilde{U}_{n}(x)$ has an $L$-extension
$\tilde{\tilde{V}}_n(x)\subset\tilde{\tilde{U}}_{n}(x)$ with
 a corresponding map $\displaystyle q_{n,x}\colon\tilde{\tilde{U}}_n(x)\to\tilde{U}_{n}(x)$  such that
$\displaystyle\big(q_{n,x}\big)|\tilde{\tilde{V}}_n(x)$ is an $L$-soft surjection onto $\tilde{V}_{n}(x)$.

Consider the nerve $\Re$ of $\nu$ and a barycentric map
$\theta\colon\Omega\to |\Re|$.  Any simplex
$\sigma=<W_0,W_1,..,W_k>$ from $\Re$, where $W_i\in\nu_{n(i)}$,
can be ordered such that $n(0)<n(1)<...<n(k)$. This is possible
because $\cap\{W_i:i=0,1,..,k\}\neq\emptyset$, so the numbers
$n(i)$ are different. It is easily seen that, for fixed $k\geq 1$
and $W\in\nu_k$,  condition (4) from the auxiliary construction
implies the following one

\begin{itemize}
\item[(5)]$St\big(\overline{g}(W\cap B), \alpha_{k}\big)\subset
V_{k-1}\big(x_W\big)$, and therefore
$St\big(f^{-1}\big(\overline{g}(W\cap B)\big),
\tilde{\alpha}_{k}\big)\subset\tilde{V}_{k-1}\big(x_W\big)$.
\end{itemize}

Let $\Sigma(\sigma)$, $\sigma\in\Re$, be the closed subset
$\theta^{-1}(\sigma)$ of $\Omega$ and
$\Sigma^k=\theta^{-1}(\Re^k)$, where  $\Re^k$ denotes the $k$-th
skeleton of $\Re$. For every  $k\geq 0$ and
$\sigma=<W_0,W_1,...,W_k>\in\Re^k$ with $W_0\in\nu_{n(0)}$, we
define by induction maps $h_k\colon \Sigma^k\to M$ and
$\displaystyle h_{\sigma}\colon\Sigma(\sigma)\to
\tilde{\tilde{U}}_{n(0)-1}(x_{W_0})$ such that

\begin{itemize}
\item[(6)] $\displaystyle h_k|\Sigma^{k-1}=h_{k-1}$ for $k\geq 1$
and $\displaystyle
h_k|\Sigma(\sigma)=q_{n(0)-1,x_{W_0}}\circ\big(h_{\sigma}|\Sigma(\sigma)\big)$
for $k\geq 0$
\end{itemize}

\noindent
and

\begin{itemize}
\item[(7)] $\displaystyle f^{-1}\big(\overline{g}(W_0\cap
B)\big)\bigcup
h_k\big(\Sigma(\sigma)\big)\subset\tilde{U}_{n(0)-1}(x_{W_0})$,
$k\geq 0$.
\end{itemize}

We also require that
\begin{itemize}
\item[(8)] $\displaystyle
h_{\sigma_1}|\big(\Sigma(\sigma_1)\cap\Sigma(\sigma_2)\big)=h_{\sigma_2}|\big(\Sigma(\sigma_1)\cap\Sigma(\sigma_2)\big)$
for any $\sigma_1$ and $\sigma_2$ from $\Re^k$ having the same
first vertex.
\end{itemize}

For $k=0$ we define  $h_0\colon \Sigma^0\to M$ and $\displaystyle
h_{<W>}\colon\Sigma(<W>)\to \tilde{\tilde{U}}_{n-1}(x_{W})$ by
$h_0\big(\Sigma(<W>)\big)=P(x_W)$ and
$h_{<W>}\big(\Sigma(<W>)\big)=Q(x_W)$, where $W\in\nu_n$ and
$Q(x_W)$ is a point from
 $\tilde{\tilde{V}}_{n-1}(x_{W})$ with $q_{0,x_W}\big(Q(x_W)\big)=P(x_W)$.
Obviously, $h_0$ restricted on every set $W\cap\Sigma^0$ is
constant, so it is continuous. Moreover, every $\displaystyle
h_{<W>}$ is also constant satisfying condition (6), and, by (5),
$h_0$ satisfies also (7).  Note that condition (8) holds for
$k=0$.

Suppose that for some $k\geq 1$ maps $h_{k-1}\colon
\Sigma^{k-1}\to M$ and $h_{\sigma}\colon\Sigma(\sigma)\to
\tilde{\tilde{U}}_{m-1}(x_{W})$ satisfying conditions (6),  (7)
and (8) have already been defined. Here $\sigma\in\Re^{k-1}$ and
$W\in\nu_m$ is the first vertex of the simplex $\sigma$.

Now, let $\sigma=<W_0,W_1,..,W_k>\in\Re^k$ with
$W_i\in\nu_{n(i)}$, $i=0,1,..,k$. Then $\sigma\cap\Re^{k-1}$
consists of the simplexes
$\sigma_i=<W_0,..,W_{i-1},W_{i+1},..,W_k>$, $i=1,2,..,k$ and the
simplex $\sigma_0=<W_1,W_2,..,W_k>$.

{\em Claim. $f^{-1}\big(\overline{g}(W_0\cap B)\big)\bigcup
h_{k-1}\big(\Sigma(\sigma_0)\big)\subset\tilde{V}_{n(0)-1}(x_{W_0})$
and $f^{-1}\big(\overline{g}(W_0\cap B)\big)\bigcup\\
h_{k-1}\big(\Sigma(\sigma_i)\big)\subset\tilde{U}_{n(0)-1}(x_{W_0})$
for every $i=1,..,k$}.

Indeed, by (7) we have $f^{-1}\big(\overline{g}(W_1\cap
B)\big)\bigcup h_{k-1}\big(\Sigma
(\sigma_0)\big)\subset\tilde{U}_{n(1)-1}(x_{W_1})$. But
$\overline{g}(W_1\cap B)\cap \overline{g}(W_0\cap
B)\neq\emptyset$, and hence $f^{-1}\big(\overline{g}(W_1\cap
B)\big)\bigcup h_{k-1}\big(\Sigma (\sigma_0)\big)$ is contained in
$St\big(f^{-1}\big(\overline{g}(W_0\cap B),
\tilde{\alpha}_{n(1)-1}\big)\big)$. Since $n(0)\leq n(1)-1$,
$\tilde{\alpha}_{n(1)-1}$ refines $ \tilde{\alpha}_{n(0)}$.  This
fact and the inclusion $St\big(f^{-1}\big(\overline{g}(W_0\cap B),
\tilde{\alpha}_{n(0)}\big)\big)\subset\tilde{V}_{n(0)-1}(x_{W_0})$,
which follows from (5), complete the proof of the claim for $i=0$.
Since $W_0$ is a vertex of each $\sigma_i$, $i=1,2,..,k$, the
other inclusions from the claim  follow directly from (7).

Consider the ``boundary" $\partial\Sigma
(\sigma)=\bigcup_{i=0}^{i=k}\Sigma(\sigma_i)$ of $\Sigma(\sigma)$.
 According to the claim, $h_{k-1}\big(\partial\Sigma (\sigma)\big)\subset\tilde{U}_{n(0)-1}(x_{W_0})$ and
$h_{k-1}\big(\overline{\partial\Sigma
(\sigma)\backslash\Sigma_0}\big)\subset\tilde{V}_{n(0)-1}(x_{W_0})$,
where $\Sigma_0=\bigcup_{i=1}^{i=k}\Sigma(\sigma_i)$. Since the
maps $\displaystyle h_{\sigma_i}\colon\Sigma(\sigma_i)\to
\tilde{\tilde{U}}_{n(0)-1}(x_{W_0})$, $i=1,..,k$, satisfy
condition (8), they determine a map $h_{\Sigma}\colon\Sigma_0\to
\tilde{\tilde{U}}_{n(0)-1}(x_{W_0})$ such that
$h_{\sigma_i}|\Sigma(\sigma_i)=h_{\Sigma}|\Sigma(\sigma_i)$ for
each $i$. Moreover, by (6), $q_{n(0)-1,x_{W_0}}\circ
h_{\Sigma}=h_{k-1}|\Sigma_0$. Therefore, we can apply Lemma 3.8
for the pair $\displaystyle
\tilde{V}_{n(0)-1}(x_{W_0})\subset\tilde{U}_{n(0)-1}(x_{W_0})$,
its $L$-extension
$\displaystyle\tilde{\tilde{V}}_{n(0)-1}(x_{W_0})\subset\tilde{\tilde{U}}_{n(0)-1}(x_{W_0})$,
the sets $\Sigma_0\subset\partial\Sigma (\sigma)\subset\Sigma
(\sigma)$ and the maps $h_{\Sigma}$ and $h_{k-1}|\partial\Sigma
(\sigma)$.  In this way we obtain a map $\displaystyle
h_{\sigma}\colon\Sigma
(\sigma)\to\tilde{\tilde{U}}_{n(0)-1}(x_{W_0})$ such that
$\displaystyle q_{n(0)-1,x_{W_0}}\circ h_{\sigma}|\partial\Sigma
(\sigma)=h_{k-1}|\partial\Sigma (\sigma)$.  Now we define
$h_k\colon \Sigma^k\to M$ by $\displaystyle h_k|\Sigma
(\sigma)=q_{n(0)-1,x_{W_0}}\circ h_{\sigma}$. Obviously, $h_k$ is
continuous on every ``simplex"  $\Sigma (\sigma)$, $\sigma\in
\Re^k$, and, since the family $\nu$ is locally finite in $\Omega$,
$h_k$ is continuous.  Moreover, $h_k$  and $h_{\sigma}$ satisfy
conditions (6), (7) and (8), and  the induction is completed.

Finally, we define $h\colon\Omega\to M$ letting $h|\Sigma^k=h_k$
for each $k$. Continuity of $h$ follows from continuity of each
$h_k$ and the fact that $\nu$ is locally finite.  Observe  also
that $(f\circ h)|p^{-1}(A)$ is $\omega$-close to $g\circ p$
because of condition (7).
\end{proof}

\begin{pro}
Let $L$ be a complex  $($not necessarily quasi-finite$)$ with the
soft map property and $f_0\colon M\to X$ be a closed map such that
each fiber $f_0^{-1}(x)$, $x\in X$, is $UV(L)$-connected in $M$.
Then  for every map $g_0\colon A\to X$, where $A$ is a closed
subset of a space $Z$ with $\e Z\leq L$ such that  either $A$  or
$g_0(A)$ is a $C$-space, there exists  a  neighborhood $Q$ of $A$
in $Z$ and an u.s.c map $\Psi\colon Q\to M $ such that $\Psi$ is
single-valued on $Q\backslash A$ and $f_0\circ\Psi$ is a
continuous single-valued map extending $g_0$.
\end{pro}

\begin{proof}
Our proof is based on some ideas from \cite[proof of Theorem
3.1]{ar}. Let  $f_0$ and $g_0$  be as in the proposition. We take
sequences $\{\omega_n\}\subset cov(X)$ and $\{\gamma_n\}\subset
cov(A)$, and open intervals $\{\bigtriangleup_n\}$ covering the
interval $J=[0,1)$, with $0\in\bigtriangleup_1$, such that:

\begin{itemize}
\item $\omega _{n+1}$ is a strong star-refinement of $\omega _n$
and $\gamma _{n+1}$ is a strong star-refinement of $\gamma _n$,
$n=1,2,3,\dots$

\item $\lim \mesh(\displaystyle\omega_n)=\lim
\mesh(\displaystyle\gamma_n)=0$

\item $\bigtriangleup_n\cap\bigtriangleup_m\neq\emptyset$ if and
only if $n$ and $m$ are consecutive integers.
\end{itemize}

\noindent Then $\omega=\{\omega_n\times\bigtriangleup_n:
n=1,2,...\}$ and $\gamma=\{\gamma_n\times\bigtriangleup_n:
n=1,2,...\}$ are open covers, respectively, of $X\times J$ and
$A\times J$, satisfying the following conditions:

\begin{itemize}
\item[$(9_i)$] For every point $(x,1)\in X\times I$ and its
neighborhood $U$ in $X\times I$ there exists another neighborhood
$V$ such that $St(V,\omega)\subset U$.
\end{itemize}

\begin{itemize}
\item[$(9_{ii})$] For every point $(a,1)\in A\times I$ and its
neighborhood $U$ in $A\times I$ there exists another neighborhood
$V$ such that $St(V,\gamma)\subset U$.
\end{itemize}

Since $f_0$ is a closed map all fibers of which are
$UV(L)$-connected in $M$, the map $f=f_0\times id\colon M\times J
\to X\times J$ has the property $(UV)$ from Lemma 3.10. Further,
let $g$ denote the map $g_0\times id\colon A\times J\to X\times J$
and consider an $L$-soft surjection $p\colon Y\to Z\times I$,
$I=[0,1]$, such that $Y$ is a space of $\e Y\leq L$.  We have the
following diagram:

$$
\begin{diagram}
\node{Y} \arrow{s,l}{\mbox{$p$ ($L$-soft)}}  \node[3]{M\times J} \arrow{s,r}{\displaystyle f=f_0\times id}\\
\node{\hspace{1.5cm}Z\times I\supset A\times J}
\arrow[3]{e,t}{\displaystyle g=g_0\times id}
 \node[3]{X\times J}
\end{diagram}
$$

Since the product of any metrizable $C$-space and $J$ is also a
$C$-space,  either $A\times J$ or $g_0(A)\times J$ is a $C$-space.
Following the notations from Lemma 3.10, we can apply construction
of this lemma by considering the spaces   $M\times J$, $X\times
J$, $Z\times J$, $A\times J$ and $p^{-1}(Z\times J)$ instead of
the spaces $M$, $X$, $Z$, $A$ and $Y$, respectively. Let us also
note that in our situation we take $\alpha_n$ and $\beta_n$,
$n\geq 0$, to be open covers of $X\times J$ satisfying condition
(1) from the auxiliary construction with $\alpha_0$ refining
$\omega$. We also require $\mu_n^*$ to be disjoint open families
in $A\times J$ satisfying condition (2) such that
$\mu^*=\bigcup_{n=1}^{\infty}\mu_n^*$ is a locally finite open
cover of $A\times J$ which, in addition, refines $\gamma$. Then,
as in the auxiliary construction, we can extend $\mu_n^*$ to
disjoint open families $\mu_n$ in $Z\times J$ by choosing an
u.s.c. retraction $r\colon Z\times I\to A\times I$ such that
$r(z,t)\subset A\times\{t\}$ for every $t\in I$. This can be
achieved by taking an u.s.c. finite-valued retraction $r_1\colon
Z\to A$ and letting $r(z,t)=r_1(z)\times\{t\}$.  Observe that this
special choice of $r$ implies that $r^{\sharp}(T)$ is open in
$Z\times I$ for every open $T\subset A\times I$ and
$r^{\sharp}(T)$ is contained in $Z\times J$ provided $T\subset
A\times J$. We also pick the points $x_W\in X\times J$, $W\in\mu$,
satisfying condition (4).

According to Lemma 3.10, there exists a map $h\colon p^{-1}(G)\to
M\times J$, where $G=\cup\{\Lambda: \Lambda\in\mu\}$, such that
each $h_k=h|\Sigma^k$ satisfies condition (7) and $(f\circ
h)|\big(p^{-1}(A\times J)\big)$ is $\omega$-close to $g\circ p$.
Now, let $H=p^{-1}\big(G\cup (A\times\{1\})\big)$ and define the
set-valued map $\psi\colon H\to M\times I$ letting $\psi(y)=h(y)$
if $y\in p^{-1}(G)$ and $\psi(y)=\big(f_0^{-1}(g_0(p(y))),1\big)$
if $y\in p^{-1}(A\times\{1\})$. Let also $\psi_1=\pi\circ\psi
\colon H\to M$, where $\pi\colon M\times I\to M$ is the
projection.

{\em Claim. The map $\psi_1$ is u.s.c.}

Since $\pi$ is continuous, it suffices to prove that $\psi$ is
u.s.c.  To this end, observe that $p^{-1}(G)$ is open in $H$ and
$\psi$ is single-valued and continuous on $p^{-1}(G)$, so that we
need to show only that $\psi$ is u.s.c. at the points of
$p^{-1}(A\times\{1\})$. Let $\{y_i\}\subset H$ be a sequence
converging to a point $y_0\in p^{-1}(A\times\{1\})$ and
$U_0=V_0\times (t,1]$ be a neighborhood of
$\psi(y_0)=\big(f_0^{-1}(g_0(p(y_0))),1\big)$ in $M\times I$.  We
are going to show that $\psi(y_i)\subset U_0$ for almost all $i$
which will complete the proof of the claim. Since $f_0$ is a
closed map, $\psi$ is u.s.c. on $p^{-1}(A\times\{1\})$. Therefore
we can assume that $\{y_i\}\subset p^{-1}(G)$, hence
$\psi(y_i)=h(y_i)$ for all $i$. Thus $p(y_0)=(a,1)\in
A\times\{1\}$ and $p(y_i)\in G$. Since $f_0$ is closed, we can
find a neighborhood $V$ of  $g_0(p(y_0))$ in $X$ with
$f_0^{-1}(V)\subset V_0$.
 By ($9_i$), there exists a neighborhood $U_1=V_1\times (q,1]$ of $\big(g_0(p(y_0)),1\big)$ in $X\times I$
such that $St(U_1,\omega)\subset U=V\times (t,1]$. Choose a
neighborhood $T(a)$ of $a$ in $A$ with $g_0(T(a))\subset V_1$ and
apply ($9_{ii}$) to find a neighborhood $S=T_1(a)\times (q^*,1]$
of $(a,1)$ in $A\times I$ such that $St(S,\gamma)\subset
T(a)\times (q,1]$. Then $r^{\sharp}(S)$ is a neighborhood of
$(a,1)$ in $Z\times I$. Since $\{p(y_i)\}$ converges to $(a,1)$,
we can assume that $\{p(y_i)\}\subset r^{\sharp}(S)$. It suffices
to show that $f(h(y_i))\in U$ for all $i$. To this end, fix $i$
and $\Lambda_{0}\in\mu_{k(0)}$ containing $p(y_i)$, where
 $k(0)$ is the minimal $k$ such that  $p(y_i)$ is contained in some element of $\mu_{k}$. Then
 $\Lambda_{0}=r^{\sharp}(\Lambda^*_0)$ for some
 $\Lambda^*_0\in\mu^*$ and therefore
$p(y_i)\in r^{\sharp}(\Lambda_{0}^*)\cap r^{\sharp}(S)$.
Consequently, $S$ meets $\Lambda_{0}^*$ and let $p(y_i^{*})\in
\Lambda_{0}^*\cap S$, where $y_i^{*}\in p^{-1}(\Lambda_0^{*})$. On
the other hand, there exists $\displaystyle\Gamma\in\gamma$
containing $\Lambda^*_0$ (recall that $\mu^*$ refines $\gamma$).
Therefore, $p(y_i^{*})\in St(S,\gamma)\subset T(a)\times (q,1]$.
Since $g\big(p(y_i^{*})\big)=(g_0\times id)\big(p(y_i^{*})\big)$,
according to the choice of $T(a)\times (q,1]$ we have

\begin{itemize}
\item[(10)]$g\big(p(y_i^{*})\big)\in U_1=V_1\times (q,1]$.
\end{itemize}

\noindent Since $k(0)$ is the minimal $k$ such that $y_i$ is
contained in some $W\in\nu_{k}$, according to the definition of
the maps $h_k$ and condition (7) from Lemma 3.10, we have
$h(y_i)\in \tilde{U}_{k(0)-1}(x_{W_0})$, where
$W_0=p^{-1}(\Lambda_0)$.  The last inclusion implies $f(h(y_i))\in
U_{k(0)-1}(x_{W_0})$.  Also, condition (5) from Lemma 3.10 yields
that

\begin{itemize}
\item[(11)] $g\big(p(y_i^*)\big)\in g\big(p\big(W_0\cap p^{-1}(A\times J)\big)\big)\subset V_{k(0)-1}(x_{W_0})$.
\end{itemize}

Hence, both $g\big(p(y_i^*)\big)$ and $f(h(y_i))$ are points from
$U_{k(0)-1}(x_{W_0})$. But the cover $\alpha_{k(0)-1}$ refines
$\omega$, and hence $U_{k(0)-1}(x_{W_0})$ is contained in an
element $O$ of $\omega$. Therefore, $O$ contains
$g\big(p(y_i^*)\big)$ and $f(h(y_i))$. This means, according to
(10), that $f(h(y_i))\in St(U_1,\omega)$. Finally, since
$St(U_1,\omega)\subset U$, we obtain $f(h(y_i))\in U$ which
completes the proof of the claim.

Now we can finish the proof. There exists a decreasing sequence
$\{Q_i\}$ of open subsets of $Z$ and an increasing sequence of
real numbers $0=t_0<t_1<..<1$ such that
$\bigcap_{i=1}^{\infty}Q_i=A$, $\lim t_i=1$,
$\overline{Q}_{i+1}\subset Q_i$ and $Q_i\times [0,t_i]\subset G$
for all $i$. Let $\varphi_i\colon Z\to [t_{i-1},t_{i}]$, $i\geq
1$, be continuous functions such that $\varphi_i(Z\backslash
Q_i)=t_{i-1}$ and $\varphi_i(z)=t_i$ for $z\in\overline{Q}_{i+1}$.
Then $\varphi\colon Z\to [0,1]$ defined by $\varphi
(z)=\varphi_i(z)$ for $z\in Q_i\backslash Q_{i+1}$, $\varphi
(Z\backslash Q_1)=0$, and $\varphi (A)=1$, is continuous.
Consequently, the map $\theta\colon Q_1\to
G\cup\big(A\times\{1\}\big)$, $\theta (z)=(z,\varphi(z))$, is well
defined and continuous. Moreover, $\theta (z)=(z,1)$ for all $z\in
A$. Since $p$ is $L$-invertible and $\e Q_1\leq L$ (as an open
subset of $Z$), we can lift $\theta$ to a map
$\overline{\theta}\colon Q_1\to H$. Then
$\Psi=\psi_1\circ\overline{\theta}\colon Q\to M$, where $Q=Q_1$,`
is the required map.
\end{proof}

Theorem 3.12 below is a generalization of the well known result
that if $G$ is an u.s.c. decomposition of a metrizable space $X$
such that each element of $G$ is $\displaystyle UV^n$ in $X$, then
$X/G$ is $LC^{n}$  \cite[Theorem 11]{da}. The result from Theorem
3.12 was also established in  \cite[Corollary 7.5]{bchk} for
finite complexes $L$ and proper $UV(L)$-maps between Polish spaces
($UV(L)$-maps are maps with all fibers being $UV(L)$-spaces). The
version of Theorem 3.12 when $L$ is a point is a generalization of
the well known result of Ancel \cite[Theorem C.5.9]{a}. This
version was also  established in \cite[Proposition 3.5]{chv1}.

\begin{thm}
Let $L$ be quasi-finite and $f\colon X\to Y$ be a closed map with
all fibers being $UV(L)$-connected in $X$. Then $Y$ is an $ANE(L)$
with respect to $C$-spaces. If, in addition, $X$ is $C^L$ $($i.e.,
every map into $X$ is $L$-homotopic to a constant map in $X$$)$,
then $Y\in AE(L)$ with respect to $C$-spaces.
\end{thm}

\begin{proof}
Let $g\colon A\to Y$ be an arbitrary map, where $A$  is a closed
subspace of a space $Z$ with $\e Z\leq L$, such that $A$ is a
$C$-space. Since $L$ is quasi-finite, it has the soft mapping
property. Therefore we can apply Proposition 3.11 to obtain a
neighborhood $U$ of $A$ in $Z$ and an u.s.c. map $\Psi\colon U\to
X$ such that $\Psi$ is single-valued outside $A$ and $f\circ\Psi$
is a single-valued extension of $g$. Hence, $Y\in ANE(L)$ with
respect to $C$-spaces (actually we proved that $Y\in ANE(g,A,Z)$
for arbitrary $g\colon A\to Y$, where $A$ is a closed subspace of
$Z$ such that $\e Z\leq L$ and $A$ is a $C$-space).

Suppose now that $X$ is $C^L$ and let $A\subset Z$ and $g\colon
A\to Y$ be as above. To show that $Y\in AE(L)$ with respect to
$C$-spaces, we need to extend  $g$ over $Z$.  Embedding $Z$ as a
closed subset of an $AE(L)$-space with $\e \leq L$, we can assume
that $Z\in AE(L)$.  Then, as before, there exists a neighborhood
$U$ of $A$ in $Z$ and an u.s.c. map $\Psi\colon U\to X$ such that
$\Psi$ is single-valued outside $A$ and $f\circ\Psi$ extends $g$.
Take neighborhoods $V_1$ and $V_2$ of $A$ in $Z$  such that
$\overline{V_1}\subset V_2\subset\overline{V_2}\subset U$.  Let
$W=Z\backslash\overline{V_1}$ and  $F=W\cap\overline{V_2}$. Since
$W\cap U$ is open in the $AE(L)$-space $Z$, the cone $\cone(W\cap
U)$ is an $AE(L)$.  So, the inclusion $F\subset W\cap U$ can be
extended to a map $\varphi\colon W\to\cone(W\cap U)$ because $F$
is closed in $W$ and $\e W\leq L$.   On the other hand, since
$X\in C^L$, $\Psi|\big(W\cap U\big)$ is $L$-homotopic to a
constant map in $X$. Consequently, the map $\Psi|F$ can be
extended to a map $h\colon W\to X$.  Finally, we define the
set-valued map $\theta\colon Z\to X$ by $\theta (z)=h(z)$ if $z\in
Z\backslash V_2$ and $\theta (z)=\Psi(z)$ otherwise. Obviously,
$\theta$ is u.s.c. and single-valued outside $A$. Moreover,
$f\circ\theta$ is the required extension of $g$.
\end{proof}

We say that a space $X$ is locally $ANE(L)$ if every point from
$X$ is $UV(L)$ in $X$. Let us mention the following corollary from
Theorem 3.12.

\begin{cor}
Let $Y$ be locally $ANE(L)$, where $L$ is quasi-finite. Then $Y\in
ANE(L)$ with respect to $C$-spaces. If, in addition, $Y\in C^L$,
then $Y\in AE(L)$ with respect to $C$-spaces.
\end{cor}

\begin{remark} We can show that if, in Corollary 3.13, the property of $X$
to be locally $ANE(L)$ is replaced by the weaker one $X$ to be
$LC^L$ (every $x\in X$ is $UV(L)$-homotopic in $X$ \cite{ch1}),
then $X$ is an $ANE(L)$ with respect to finite-dimensional spaces
(see also \cite[Theorem 4.1]{bchk} for a similar result).
\end{remark}

We know that the domain and the range of a $UV^n$-map between
compacta are simultaneously $UV^n$  (see, for example \cite{be}).
Here is a generalization of this result for a subclass of
quasi-finite complexes.  We  say that a $CW$ complex $L$ is a
$C$-complex if every space of $\e \leq L$ is a $C$-space.  Each
complex $L$ with $L\leq\bbb S^n$  for some $n$ (this means that
$\e  Z\leq L$ implies $\dim Z\leq n$ for any space $Z$) is a
$C$-complex, in particular every sphere $\bbb S^k$ is such a
complex.  Observe that Lemma 3.10 and Proposition 3.11 remain
valid for $C$-complexes $L$ having the soft map property without
the requirements  either $A$ or $g(A)$ (resp., $g_0(A)$) to be
$C$-spaces.  This yields that, if in Theorem 3.12 and Corollary
3.13 $L$ is a quasi-finite $C$-complex, then $Y$ is an $A(N)E(L)$.

\begin{thm}
Let $L$ be a quasi-finite $C$-complex and $f\colon X\to Y$ a
closed map with $UV(L)$-fibers. Then $X$ is $UV(L)$ if and only if
$Y$ is.
\end{thm}

\begin{proof}
Let $E_X$ be a normed space containing $X$ as a strong $Z$-set.
This means that $X\subset E_X$ is closed and for every $\omega\in
cov(E_X)$ and every map $g\colon Q\to E_X$, where $Q$ is an
arbitrary space, there is another map $h\colon Q\to E_X$ which is
$\omega$-close to $g$ and $\overline{h(Q)}\cap X=\emptyset$ (such
space $E_X$ can be constructed as follows: embed $X$ as a closed
subset of a normed space $F$ and let $E_X$ be the product $F\times
l_2(\tau)$, where $w(X)\leq\tau$; then $X\times\{0\}$ is a copy of
$X$ which is a strong $Z$-set in $E_X$). Identifying each fiber of
$f$ with a point, we obtain space $E_Y$ (equipped with the
quotient topology) and let $p\colon E_X\to E_Y$ be the natural
quotient map. Obviously, $p(X)\subset E_Y$ is closed and, since
$f$ is a closed map, $p(X)$ is homeomorphic to $Y$. And everywhere
below we write $Y$  instead of $p(X)$. Moreover, $p$ is a closed
map and $E_Y$ is metrizable.  Any fiber of $p$ is either a point
or $f^{-1}(y)$ for some $y\in Y$. Hence, $p$ is an $UV(L)$-map.
Since $E_X$ is an absolute extensor for metrizable spaces, the
fibers of $p$ are $UV(L)$-connected in $E_X$. Consequently, by the
modified version of Theorem 3.12 for $C$-complexes, $E_Y\in
AE(L)$.

$X\in UV(L)\Rightarrow Y\in UV(L)$. To prove this implication, by
Corollary 3.7, it suffices to show that $Y$ is $UV(L)$ in $E_Y$.
Let $U$ be a neighborhood of $Y$ in $E_Y$.  Since $X$ is $UV(L)$
in $E_X$ (recall that $E_X$ is an absolute extensor) and $p$ is
closed, there exists a neighborhood $V$ of $Y$ in $E_Y$ such that
the pair $p^{-1}(V)\subset p^{-1}(U)$ is $L$-connected. We choose
a neighborhood $V_1$ of $Y$ in $E_Y$ with $\overline{V}_1\subset
V$ and  show that the pair $V_1\subset U$ is  $L$-connected. To
this end,  take a space $Z$ with $\e Z\leq L$ and a map $h\colon
A\to V_1$ with $A\subset Z$ being closed. Since $U$ is an
$ANE(L)$, there exists $\omega\in cov (U)$ satisfying condition
(H) from Proposition 3.2.  Further,  let $\beta\in cov(E_Y)$ be
the cover $\{G\cap V:
G\in\omega\}\cup\{E_Y\backslash\overline{V}_1\}$. By Lemma 3.10,
there exists a map $h_1\colon A\to E_X$ such that $p\circ h_1$ is
$\beta$-close to $h$. Obviously, $h_1(A)\subset p^{-1}(V)$ and
hence there exists an extension $h_2\colon Z\to p^{-1}(U)$ of
$h_1$. Then $p\circ h_2$ is a map from $Z$ into $U$ such that
$(p\circ h_2)|A$ is $\omega$-close to $h$. Finally, according to
the choice of $\omega$, $h$ admits an extension
$\overline{h}\colon Z\to U$.

$Y\in UV(L)\Rightarrow X\in UV(L)$. As in the previous
implication, it suffices to show that $X$ is $UV(L)$ in $E_X$.  To
this end, let $U$ be a neighborhood of $X$ in $E_X$.  We can
assume that $U=p^{-1}(U_0)$ for some neighborhood $U_0$ of $Y$ in
$E_Y$. Choose neighborhoods $V_0$, $G_0$ and $W_0$ of $Y$ such
that $V_0\subset\overline{V_0}\subset
G_0\subset\overline{G_0}\subset W_0\subset\overline{W_0}\subset
U_0$ and the pair $G_0\subset W_0$ is $L$-connected. Denote by
$V$, $G$ and $W$, respectively,  the preimages $p^{-1}(V_0)$,
$p^{-1}(G_0)$ and $p^{-1}(W_0)$. We claim that  the pair $V\subset
U$ is  $L$-connected.  Indeed, consider a map $g_V\colon A\to V$,
where $A$ is a closed subset of a space $Z$ with $\e Z\leq L$. Let
$\alpha\in cov (U)$ satisfy condition (H) from Proposition 3.2 and
$\alpha_1=\{T\cap G:
T\in\alpha\}\cup\{E_X\backslash\overline{V}\}\in cov(E_X) $. Since
$X$ is a strong $Z$-set in $E_X$, we can find a map $g_G\colon
A\to E_X$  which is $\alpha_1$-close to $g_V$ and
$\overline{g_G(A)}\cap X=\emptyset$.  It is easily seen that
$g_G(A)\subset G$ and $g_G$  is $\alpha$-close to $g_V$.  The last
yields (because of the choice of  $\alpha$) that   $g_V$ can be
extended to a map from $Z$ into $U$ if and only if $g_G$ has such
an extension. Hence, our proof is reduced to show that $g_G$
admits an extension from $Z$ into $U$.  Obviously, $g_G$ can be
considered as a map from $A$ into $G_0$ such that the closure
$\overline{g_G(A)}$  (this is a closure in $E_Y$) does not meet
$Y$. Since $G_0\subset W_0$ is $L$-connected, $g_G$ can be
extended to a map $g_W\colon Z\to W_0$. Finally, consider the
cover $\gamma\in cov(E_Y)$ defined by $\gamma=\{p(T\backslash X):
T\in\alpha\}\cup\{E_Y\backslash\overline{g_G(A)}\}\cup\{E_Y\backslash\overline{W_0}\}$.
According to Lemma 3.10, there exists a map $g_U\colon Z\to E_X$
such that $p\circ g_U$ is $\gamma$-close to $g_W$.  It is easily
seen that $g_U(Z)\subset U$ and $g_U|A$ is $\alpha$-close to
$g_G$. The last condition implies that $g_G$ admits an extension
from $Z$ into $U$ which completes our proof.
\end{proof}


\bigskip

\end{document}